\documentclass[a4paper,11pt,twoside,reqno]{amsart}

\usepackage[utf8]{inputenc}
\usepackage[plainpages=false,pdfpagelabels=true]{hyperref}
\usepackage{amssymb,amsthm}
\usepackage[margin=1in]{geometry}

\newtheorem{Satz}{Theorem}[section]
\newtheorem{Prop}[Satz]{Proposition}
\newtheorem{Lem}[Satz]{Lemma}

\newtheorem{Cor}[Satz]{Corollary}

\newcommand{\vol}{{\operatorname{Vol}}}
\theoremstyle{definition}
\newtheorem{Dfn}[Satz]{Definition}
\newtheorem{Bem}[Satz]{Remark}

\newcommand{\cM}{{\mathcal{M}}}

\allowdisplaybreaks[1]

\renewcommand{\epsilon}{\varepsilon}

\numberwithin{equation}{section}

\title{The normalized second order renormalization group flow on closed surfaces}
\author{Volker Branding}
\date{\today}
\address{TU Wien\\
Institut für diskrete Mathematik und Geometrie\\
Wiedner Hauptstraße 8–10, A-1040 Wien}
\email[]{volker@geometrie.tuwien.ac.at}
\begin{document}

\begin{abstract}
We study a normalized version of the second order renormalization group flow
on closed Riemannian surfaces. We discuss some general properties of this flow
and establish several basic formulas. In particular, we focus on surfaces
with zero and positive Euler characteristic.
\end{abstract} 

\maketitle

\section{Introduction and Results}
Harmonic maps from Riemannian surfaces to Riemannian manifolds are a prominent variational problem
in differential geometry. However, they also have a dual live in theoretical physics.
Namely, the non-linear bosonic sigma model is governed by the action
\begin{equation}
\label{polyakov-action}
E(\phi)=\frac{1}{4\pi\alpha'}\int_M\frac{\partial\phi^i}{\partial x_\alpha}\frac{\partial\phi^j}{\partial x_\beta}g_{ij}h^{\alpha\beta}dM.
\end{equation}
Here, \((M,h_{\alpha\beta})\) and \((N,g^{ij})\) are Riemannian manifolds and \(\phi\colon M\to N\) is a map.
Moreover, the Einstein summation convention is applied.
The constant \(\alpha'>0\) is usually interpreted as coupling constant.

The perturbative quantization of the action \eqref{polyakov-action} requires to introduce a 
momentum cutoff \(\Lambda>0\) and one is left with a family of quantum field theories
depending on \(\Lambda\). The requirement that this family should be independent
of the cutoff leads to the renormalization group equation:

\[
\frac{\partial g_{ij}}{\partial\Lambda}=-\beta_{ij}.
\]
The function \(\beta\) is determined by a perturbative expansion:
\[
\beta_{ij}=\alpha'\beta_{ij}^{(1)}+\alpha'^2\beta_{ij}^{(2)}+\alpha'^3\beta_{ij}^{(3)}+\ldots
\]
For the non-linear bosonic sigma model \eqref{polyakov-action} the beta function can be computed as
\begin{align}
\label{beta-function}
\beta_{ij}^{(1)}=&R_{ij}, \\
\nonumber \beta_{ij}^{(2)}=&\frac{1}{2}R_{iklm}R_{j}^{~klm},\\
\nonumber\beta_{ij}^{(3)}=&\frac{1}{8}\nabla_p R_{iklm}\nabla_p R_{j}^{~klm}-\frac{1}{16}\nabla_i R_{klmp}\nabla_jR^{klmp} 
+\frac{1}{2}R_{klmp}R_i^{~mlr}R_{j\hspace{0.3cm}r}^{~kp}-\frac{3}{8}R_{iklj}R^{kspr}R^l_{~spr}.
\end{align}
Here, \(R_{iklm}\) denote the components of the curvature tensor and \(R_{ij}\) represent
the components of the Ricci tensor.
For a derivation of the one and two loop order see \cite{Friedan:1980jf} and \cite{Friedan:1980jm}.
The third loop order was considered in \cite{Graham1987543} and the fourth loop order in \cite{MR1009743}.
The one loop order can also be derived in a rigorous mathematical fashion, see \cite{1408.4466}.
For a general introduction to the theory of renormalization see the book \cite{MR2778558}.

It is well-known that the first order contribution in \eqref{beta-function} gives rise to the Ricci flow
equation, which by now is probably the most famous geometric evolution equation.
However, one can of course also include the higher order terms from \eqref{beta-function} and study the resulting geometric flow.
Considering also the second order contribution from \eqref{beta-function},
one obtains the following evolution equation for the metric \(g_{ij}\):
\begin{equation}
\label{second-order-evolution}
\frac{\partial}{\partial t}g_{ij}=-\alpha'R_{ij}-\frac{\alpha'^2}{2}R_{iklm}R_{j}^{~klm}.
\end{equation}
We call \eqref{second-order-evolution} \emph{the second order renormalization group flow}.
This evolution equation has already been analyzed from a geometric perspective:
On closed surfaces the evolution equation \eqref{second-order-evolution} was studied in \cite{MR2532057}.
The case of three-dimensional homogeneous geometries is studied in \cite{MR3043753} (see also \cite{MR2822454}) and 
the authors point out that in some cases the qualitative behaviour of the second
order renormalization group flow differs from the Ricci flow.
A stability analysis using techniques of maximal regularity was performed in \cite{MR2415546}.
For a general introduction to this subject see \cite{MR3202021}.

We will study \eqref{second-order-evolution} on a closed Riemannian surface,
hence we apply the following relations between
the curvature tensor \(R_{ijkl}\), the Ricci tensor \(R_{ij}\)
and the scalar curvature \(R\)
\begin{align*}
R_{ijkl}=\frac{1}{2}R(g_{ik}g_{jl}-g_{il}g_{jk}),\qquad R_{ij}=\frac{1}{2}Rg_{ij}.
\end{align*}
In addition, we rescale the time parameter by \(t\to\frac{\alpha'}{2}t\) and thus \eqref{second-order-evolution} simplifies to
\begin{equation}
\label{2nd-rg}
\frac{\partial}{\partial t}g_{ij}=-(R+\frac{\alpha'}{4}R^2)g_{ij}. 
\end{equation}

The second order normalization group flow \eqref{2nd-rg} does not preserve the volume of the surface \(M\).
It is the aim of this article to study the following volume-preserving variant:
\begin{align}
\label{second-order-flow}
\frac{\partial}{\partial t}g=-(R+\frac{\alpha'}{4}R^2-r)g,\qquad g(0)=g_0.
\end{align}
Here, the normalization \(r\) is given by
\begin{align}
r:=\frac{\int_M (R+\frac{\alpha'}{4}R^2)d\mu}{\int_Md\mu},
\end{align}
where \(d\mu\) denotes the surface measure. The investigation of this flow was already proposed in \cite{MR3202021}.
We call \eqref{second-order-flow} the \emph{normalized second order renormalization group flow}.
For the study of this flow we consider the two subspaces of metrics
\begin{align*}
\cM_+:=\{g\in M\mid 1+\alpha'R/2>0\}, \qquad 
\cM_-:=\{g\in M\mid 1+\alpha'R/2<0\}.
\end{align*}
In our analysis of \eqref{second-order-flow}, we apply the tools for the normalized Ricci flow on closed surfaces,
see \cite{MR954419} and \cite{MR2061425}. However, due to the non-linear nature of \eqref{second-order-flow} most of them can no longer be utilized.
In the limit \(\alpha'\to 0\) the evolution equation \eqref{second-order-flow}
reduces to the normalized Ricci flow equation. 

For \(g\in\cM_+\) the normalized second order renormalization group flow behaves similar to the normalized Ricci flow.
Thus, from a mathematical perspective it would be interesting to know if one can also prove the uniformization theorem
using the normalized second order renormalization group flow. 
However, this question cannot be answered with the tools developed in this paper.

This paper is organized as follows.
In Section 2 we establish some basic formulas for \eqref{second-order-flow}. 
Afterwards, in Section 3, we discuss different curvature potentials and
apply these in the case of surfaces with zero Euler characteristic.
In the last section we focus on surfaces with positive Euler characteristic.

\section{The normalized second order renormalization group flow}
We start this section by making some general comments on the evolution equation \eqref{second-order-flow}.
\begin{Bem}
By the Gauss-Bonnet Theorem it follows that
\begin{align*}
\vol(M,g)+\alpha'\pi\chi(M)>0 \text{ on } \cM_+, \qquad 
\vol(M,g)+\alpha'\pi\chi(M)<0 \text{ on } \cM_-.
\end{align*}
Hence, we may follow that the set \(\cM_-\) is empty for \(M\approx S^2,T^2\).	 
\end{Bem}

Since we are concerned with the study of \(\eqref{second-order-flow}\) on surfaces
with \(\chi(M)=0\) and \(\chi(M)=2\), we will mostly assume that \(g\in \cM_+\).
For \(g\in \cM_+\) the normalized second order renormalization group flow \eqref{second-order-flow}
is parabolic, whereas for \(g\in \cM_-\) it is backward parabolic.

\begin{Bem}
The fixed points of \eqref{second-order-flow} satisfy
\[
R+\frac{\alpha'}{4}R^2=r 
\]
and by differentiating we obtain:
\[
\nabla R(1+\frac{\alpha'^2}{2}R)=0
\]
Hence, as for the Ricci flow on surfaces the fixed points of \eqref{second-order-evolution} are exactly the metrics of constant curvature.
\end{Bem}

Throughout this article we will make use of the following facts:
Under a conformal change of the metric \(g=e^{2u}h\), we have the following formula relating the scalar curvatures
\begin{equation}
\label{conformal-change-R}
R_g=e^{-2u}(R_h-2\Delta_h u).
\end{equation}

In addition, under \eqref{second-order-flow} the following identities hold (see for example \cite{MR2061425}, Chapter 3 and Cor.5.5):
\begin{align}
\label{evolution-christoffel}
\frac{\partial}{\partial t}\Gamma^k_{ij}&=\frac{1}{2}g^{kl}\big[\nabla_i\big(\frac{\partial}{\partial t}g_{jl}\big)+\nabla_j\big(\frac{\partial}{\partial t}g_{il}\big)
-\nabla_l\big(\frac{\partial}{\partial t}g_{ij}\big)\big], \\
\label{evolution-laplace}
\frac{\partial}{\partial t}\Delta&=\big(R+\frac{\alpha'}{4}R^2-r\big)\Delta, \\
\label{evolution-dmu}
\frac{\partial}{\partial t}d\mu&=-\big(R+\frac{\alpha'}{4}R^2-r\big)d\mu.
\end{align}
Here, \(\Gamma^k_{ij}\) are the Christoffel symbols of the metric \(g_{ij}\) and \(d\mu\) represents the surface measure on \(M\).
Moreover, we will often make use of the Ricci identity:
\begin{equation}
\nabla\Delta=\Delta\nabla-\frac{1}{2}R\nabla. 
\end{equation}

\begin{Bem}
The normalized second order renormalization group flow is not invariant under rescaling.
Suppose we rescale the metric by a positive factor \(\lambda\), more precisely \(g\to\lambda g\),
then
\[
R\to\frac{1}{\lambda }R,\qquad d\mu\to\lambda d\mu,\qquad \vol(M,g)\to\lambda\vol(M,g)
\]
and thus
\[
R+\frac{\alpha'}{4}R^2-\frac{\int_MRd\mu}{\int_Md\mu}-\frac{\alpha'}{4}\frac{\int_M R^2d\mu}{\int_Md\mu}\to
\frac{1}{\lambda}\big(R+\frac{\alpha'}{4}\frac{R^2}{\lambda}-\frac{\int_MRd\mu}{\int_Md\mu}-\frac{\alpha'}{4}\frac{1}{\lambda}\frac{\int_M R^2d\mu}{\int_Md\mu}\big),
\]
which shows that the terms on the right hand side of \eqref{second-order-flow} do not all scale the same way.
Thus, we expect that the normalized second order renormalization group flow behaves differently than
the unnormalized version. 

Moreover, contrary to the Ricci flow, the normalized second order renormalization
group flow cannot be obtained from the second order renormalization group flow by rescaling the metric.
\end{Bem}

Let us also briefly recall the following (Lem. 5.7 in \cite{MR2061425})
\begin{Lem}
If \(g\) is a smooth 1-parameter family of metrics on a Riemannian surface \(M\)
satisfying
\[
\frac{\partial g}{\partial t}=fg
\]
then the scalar curvature evolves by
\begin{equation}
\frac{\partial R}{\partial t}=-\Delta f-Rf.
\end{equation}
\end{Lem}

Using the standard theory for parabolic partial differential equations we obtain (similar to Theorem 3.1. in \cite{MR2532057})
\begin{Prop}[Short-time existence]
There exists a unique, smooth solution of \eqref{second-order-flow} in the set \(\cM_+\)
for \(t\in[0,T)\).
\end{Prop}
\begin{proof}
By \eqref{second-order-flow} we are deforming a given metric \(g(0)=h_0\)
in its conformal class, thus we may write \(g(t)=e^{2u(t)}h_0\)
with a time independent metric \(h_0\). Using this \eqref{second-order-flow} is equivalent to
\begin{equation}
\label{second-order-flow-u}
\frac{\partial u}{\partial t}=-\frac{1}{2}\big(R+\frac{\alpha'}{4}R^2-r\big)(u),\qquad u(0)=0.
\end{equation}
With the help of \eqref{conformal-change-R} we calculate the linearization of the right hand side using \(\frac{\partial u}{\partial t}\big|_{t=0}=f\), namely
\begin{align}
\frac{\partial}{\partial t}\big|_{t=0}\big(R+\frac{\alpha'}{4}R^2-r\big)(u)=&-\big(1+\frac{\alpha'}{2}R\big)(\Delta f+Rf)\\
\nonumber&+\frac{\alpha'}{4}\frac{1}{\vol(M,g)}\int_MR^2fd\mu-\frac{\alpha'}{2}\frac{1}{\vol(M,g)}\int_M\nabla R\nabla fd\mu.
\end{align}
Thus, since \(g(t)=e^{2u(t)}h_0\in \cM_+\), equation \eqref{second-order-flow-u} is \emph{parabolic}.
The existence of a short-time solution then follows from standard theorems, 
see for example \cite{MR2744149}, Prop. 8.1, p. 389.
\end{proof}

\begin{Bem}
In the case of a higher dimensional manifold,
the short-time existence of the second order renormalization group flow
has been established under the condition \(1+\alpha'K>0\) for all sectional curvatures \(K\).
See \cite{1306.1721} for the three-dimensional case and \cite{1401.1454} for the general case.
\end{Bem}

To compute the evolution of the scalar curvature, we use that if
\[
\frac{\partial g}{\partial t}=h,
\]
then
\begin{equation}
\label{evolution-R-general}
\frac{\partial R}{\partial t}=-\Delta\operatorname{Tr}h+\operatorname{div}(\operatorname{div} h)-\langle h,\operatorname{Ric}\rangle .
\end{equation}

\begin{Lem}
Let \((M,g_t)\) be a solution of the normalized second order renormalization group flow \eqref{second-order-flow}.
Then the scalar curvature \(R\) evolves by
\begin{equation}
\label{evolution-R}
\frac{\partial R}{\partial t}=(1+\frac{\alpha'}{2}R)\Delta R+\frac{\alpha'}{2}|\nabla R|^2+R(R+\frac{\alpha'}{4}R^2-r).
\end{equation}
\end{Lem}

\begin{Bem}
The evolution of the scalar curvature \eqref{evolution-R} has some similarity with the porous medium equation
\[
\frac{\partial u}{\partial t}=\Delta u^m
\]
for \(m=2\).
\end{Bem}

\begin{Bem}
If we want to use the maximum principle to obtain an estimate on how the scalar curvature \(R\) behaves under \eqref{second-order-flow} with \(g\in \cM_+\),
we have to study the comparison ODE
\begin{equation}
\label{abel}
 y'(t)=-r(t)y(t)+y^2(t)+\frac{\alpha'}{4}y^3(t).
\end{equation}
This equation is known as \emph{Abel equation of the first kind}. 
Unfortunately, it cannot be integrated directly to obtain an estimate.

In the case of the unnormalized second order renormalization group flow the analysis
of the corresponding ODE \eqref{abel} is the cornerstone for the results presented in \cite{MR2532057}.
However, these results rely heavily on the calculation of the fixed points of \eqref{abel}.
Due to the time-dependent normalization factor \(r(t)\) we cannot apply the same method here.
\end{Bem}

\begin{Cor}
Let \((M,g_t)\) be a solution of the normalized second order renormalization group flow \eqref{second-order-flow}.
Then the square of the scalar curvature \(R\) evolves by
\begin{equation}
\label{evolution-R-squared}
\frac{\partial}{\partial t}R^2=\big(1+\frac{\alpha'}{2}R\big)\Delta R^2-2|\nabla R|^2+2R^2(R+\frac{\alpha'}{4}R^2-r).
\end{equation}
\end{Cor}

This allows us to derive a first estimate:
\begin{Lem}
Let \((M,g_t)\) be a solution of the normalized second order renormalization group flow \eqref{second-order-flow} with \(g\in \cM_+\).
Then the following estimate holds:
\begin{equation}
R^2\leq R_0^2\frac{c_1e^{c_1t}}{c_1-c_2R_0^2(e^{c_1t}-1)}\textrm{  with  } c_1=1-\frac{4\pi\chi(M)}{\vol(M,g)},\qquad c_2=1+\frac{\alpha'}{2}.
\end{equation}
\end{Lem}
\begin{proof}
Using that \(g\in \cM_+\) we can apply the maximum principle to \eqref{evolution-R-squared} and estimate
\begin{align*}
\frac{\partial}{\partial t}R^2&=2R^3+\frac{\alpha'}{2}R^4-2R^2r 
\leq (1-2\frac{2\pi\chi(M)}{\vol(M,g)})R^2+(1+\frac{\alpha'}{2})R^4
\end{align*}
and the claim follows by solving the ODE.
\end{proof}

\begin{Lem}
Let \((M,g_t)\) be a solution of the normalized second order renormalization group flow \eqref{second-order-flow}.
Then \(|\nabla R|^2\) evolves by
\begin{align}
\label{evolution-nabla-R}
\frac{\partial}{\partial t}|\nabla R|^2=&(1+\frac{\alpha'}{2}R)\Delta|\nabla R|^2
-2(1+\frac{\alpha'}{2}R)|\nabla^2R|^2+\alpha'\langle\nabla|\nabla R|^2,\nabla R\rangle\\
\nonumber &
+(4R-3r+\frac{5}{4}\alpha'R^2+\alpha'\Delta R)|\nabla R|^2.
\end{align}
\end{Lem}
\begin{proof}
By a direct calculation one finds
\begin{align*}
\frac{\partial}{\partial t}\nabla R&=
\nabla\big((1+\frac{\alpha'}{2}R)\Delta R+\frac{\alpha'}{2}|\nabla R|^2+R(R+\frac{\alpha'}{4}R^2-r)\big)\\
&=\big(1+\frac{\alpha'}{2}R\big)\nabla\Delta R+\frac{\alpha'}{2}\nabla R\Delta R
+\frac{\alpha'}{2}\nabla|\nabla R|^2+2R\nabla R+\frac{3\alpha'}{4}R^2\nabla R-r\nabla R.
\end{align*}
Moreover, we find
\begin{align*}
\frac{\partial}{\partial t}|\nabla R|^2=\frac{\partial}{\partial t}(g^{ij}\nabla_iR\nabla_jR) 
=(R+\frac{\alpha'}{4}R^2-r)|\nabla R|^2 +2\langle\frac{\partial}{\partial t}\nabla R,\nabla R\rangle
\end{align*}
and combining both equations yields the result.
\end{proof}

\begin{Bem}
The evolution equation for \(|\nabla R|^2\), \eqref{evolution-nabla-R}, shows a problem that one has to deal with
when one wants to derive estimates via the maximum principle for the second order renormalization group flow.
Namely, since the Laplacian always comes with a prefactor \(1+\frac{\alpha'}{2}R\)
one gets terms in the Bochner formulas that are to nonlinear to derive suitable estimates.
\end{Bem}

In the following we will often need the evolution of the Christoffel symbols.
From \eqref{evolution-christoffel} and \eqref{second-order-flow} we obtain
\begin{equation}
\label{evolution-christoffel-second-order}
\frac{\partial}{\partial t}\Gamma^k_{ij}=\frac{1}{2}\big(1+\frac{\alpha'}{2}R\big)\big(-\nabla_iR\delta^k_j-\nabla_jR\delta^k_i+\nabla^kRg_{ij}\big).
\end{equation}

\begin{Lem}
Let \((M,g_t)\) be a solution of the normalized second order renormalization group flow \eqref{second-order-flow}.
Then \(|\nabla^2 R|^2\) evolves by
\begin{align}
\label{evolution-nabla2-R}
\frac{\partial}{\partial t}|\nabla^2 R|^2=&\big(1+\frac{\alpha'}{2}R\big)\big(\Delta|\nabla^2 R|^2
-2|\nabla^3 R|^2+2R|\nabla^2R|^2\big)+(2R^2-4r)|\nabla^2 R|^2\\
\nonumber &
+\alpha'(\Delta R|\nabla^2R|^2+2\langle\nabla R\nabla\Delta R,\nabla^2R\rangle
+\langle\nabla^2|\nabla R|^2,\nabla^2R\rangle)\\
\nonumber&+\big(4+3\alpha'R\big)\langle\nabla R\nabla R,\nabla\nabla R\rangle.
\end{align}
\end{Lem}
\begin{proof}
We calculate
\[
\frac{\partial}{\partial t}\nabla_i\nabla_jR=\nabla_i\nabla_j\big(\frac{\partial R}{\partial t}\big)-\big(\frac{\partial}{\partial t}\Gamma^k_{ij}\big)\nabla_kR
\]
and using the evolution equation for the scalar curvature \eqref{evolution-R}, we find for the first term
\begin{align*}
\nabla_i\nabla_j\big(\frac{\partial R}{\partial t}\big)=\nabla_i\nabla_j\big(\big(1+\frac{\alpha'}{2}R\big)\Delta R\big)
+\frac{\alpha'}{2}\nabla_i\nabla_j|\nabla R|^2+\nabla_i\nabla_jR^2+\frac{\alpha'}{4}\nabla_i\nabla_jR^3-r\nabla_i\nabla_jR.
\end{align*}
To manipulate the first term on the right hand side we use
\[
\nabla_i\nabla_j\Delta R=\Delta\nabla_i\nabla_jR-2R\nabla_i\nabla_jR+\big(R\Delta R+\frac{1}{2}|\nabla R|^2\big)g_{ij}-\nabla_iR\nabla_jR
\]
and find
\begin{align*}
\nabla_i\nabla_j\big(\big(1+\frac{\alpha'}{2}R\big)\Delta R\big)=&\big(1+\frac{\alpha'}{2}R\big)(\Delta\nabla_i\nabla_jR-2R\nabla_i\nabla_jR+\big(R\Delta R+\frac{1}{2}|\nabla R|^2\big)g_{ij}-\nabla_iR\nabla_jR)\\
&+\frac{\alpha'}{2}\big(\nabla_i\nabla_jR\Delta R+\nabla_jR\nabla_i\Delta R+\nabla_iR\nabla_j\Delta R\big).
\end{align*}
By the formula for the evolution of the Christoffel symbols \eqref{evolution-christoffel-second-order} we get
\begin{align*}
\big(\frac{\partial}{\partial t}\Gamma^k_{ij}\big)\nabla_kR
=&-\big(1+\frac{\alpha'}{2}R\big)\nabla_iR\nabla_jR+\frac{1}{2}\big(1+\frac{\alpha'}{2}R\big)|\nabla R|^2g_{ij}
\end{align*}
and in addition we have
\begin{align*}
\nabla_i\nabla_jR^2=&2\nabla_iR\nabla_jR+2R\nabla_i\nabla_jR,\qquad
\frac{\alpha'}{4}\nabla_i\nabla_jR^3=\frac{3}{2}\alpha'R\nabla_iR\nabla_jR+\frac{3}{4}\alpha'R^2\nabla_i\nabla_jR.
\end{align*}

Combining the formulas we find
\begin{align*}
\frac{\partial}{\partial t}\nabla_i\nabla_jR=&\big(1+\frac{\alpha'}{2}R\big)\big(\Delta\nabla_i\nabla_jR+R\Delta Rg_{ij}\big)
-\frac{\alpha'}{4}R^2\nabla_i\nabla_jR-r\nabla_i\nabla_jR+2\nabla_iR\nabla_jR\\
&+\frac{\alpha'}{2}\big(\nabla_i\nabla_jR\Delta R+\nabla_jR\nabla_i\Delta R+\nabla_iR\nabla_j\Delta R+\nabla_i\nabla_j|\nabla R|^2\big)+\frac{3}{2}\alpha'R\nabla_iR\nabla_jR.
\end{align*}
Using
\begin{align*}
\frac{\partial}{\partial t}|\nabla^2 R|^2&=\frac{\partial}{\partial t}\big(g^{ik}g^{jl}\nabla_i\nabla_jR\nabla_k\nabla_lR\big) 
=2(R+\frac{\alpha'}{4}R^2-r)|\nabla^2R|^2+2\langle\frac{\partial}{\partial t}\nabla^2 R,\nabla^2 R\rangle
\end{align*}
then gives the result.
\end{proof}

\begin{Bem}
For \(g\in \cM_+\) we can apply the maximum principle to \eqref{evolution-nabla2-R} and then have to estimate	
\begin{align*}
\frac{\partial}{\partial t}|\nabla^2 R|^2\leq\alpha'|\nabla^2R|^3+C|\nabla^2R|^2f(R,\nabla R),
\end{align*}
where \(f(R,\nabla R)\) is a function only depending on \(R,\nabla R\). However,
due to the presence of the \(|\nabla^2R|^3\) term we cannot expect to succeed in deriving an estimate using the maximum principle.
\end{Bem}

\begin{Lem}
Let \((M,g_t)\) be a solution of the normalized second order renormalization group flow \eqref{second-order-flow}.
Then the normalization function \(r\) evolves as
\begin{equation}
\label{evolution-r}
\frac{\partial r}{\partial t}=\frac{\alpha'}{4\vol(M)}\big(-2\int_M(1+\frac{\alpha'}{2}R)|\nabla R|^2d\mu+\int_M(R^3+\frac{\alpha'}{4}R^4)d\mu-r\int_M R^2d\mu\big).
\end{equation}
\end{Lem}

\begin{Bem}
We do not get an estimate from \eqref{evolution-r}, however, 
we can estimate the normalization function \(r\) via Young's equality
\begin{align}
r\geq-\frac{1}{\alpha'},\qquad r\geq\frac{2\pi\chi(M)}{\vol(M,g)}.
\end{align}
\end{Bem}

We finish this section by commenting on the the third-loop contribution to \eqref{beta-function}.
\begin{Bem}
Considering also the third order contribution from \eqref{beta-function}, then on a Riemannian surface the metric \(g_{ij}\) evolves by
\begin{align}
\label{third-order-flow-surface}
\frac{\partial}{\partial t}g_{ij}=-\alpha'\frac{R}{2}g_{ij}-\alpha'^2\frac{R^2}{8}g_{ij}-\alpha'^3\big(\frac{1}{16}|\nabla R|^2g_{ij}-\frac{1}{16}\nabla_iR\nabla_jR+\frac{5}{32}R^3g_{ij}\big).
\end{align}
Using \eqref{evolution-R-general}, we may compute the evolution of the scalar curvature under \eqref{third-order-flow-surface} and find
\begin{align*}
\label{evolution-R-third}
\frac{\partial R}{\partial t}=&\big(\frac{\alpha'}{2}+\frac{\alpha'^2}{4}R+\alpha'^3\frac{15}{16}R^2\big)\Delta R
+\big(\frac{\alpha'^2}{4}+\alpha'^3\frac{25}{32}R\big)|\nabla R|^2\\
\nonumber&+\frac{\alpha'^3}{16}\big(\frac{1}{2}R^2\Delta R+|\nabla\nabla R|^2+\frac{1}{2}R|\nabla R|^2+2\langle\nabla R,\nabla\Delta R\rangle+|\Delta R|^2\big)\\
\nonumber&+\frac{\alpha'}{2}R^2+\frac{\alpha'^2}{8}R^3+\alpha'^3\frac{5}{32}R^4.
\end{align*}
One clearly sees that including the third order contribution again changes the type of partial differential equation.
\end{Bem}

\section{Second order gradient Ricci Solitons and Curvature Potentials}
In the context of the second order renormalization group flow Ricci solitons have already been considered in \cite{MR3202021}, Section 4.
We call a solution of the second order normalized renormalization group flow \eqref{second-order-flow} \emph{self-similar}
if there exists a one-parameter family of conformal diffeomorphisms \(\varphi(t)\) such that
\begin{equation}
g(t)=\varphi(t)^*g_0.
\end{equation}
Differentiating with respect to \(t\) yields
\[
\frac{\partial}{\partial t}g=\mathcal{L}_Xg,
\]
where \(X(t)\) is the one-parameter family of vector fields generated by \(\varphi(t)\).
Using \eqref{second-order-flow} this yields
\begin{equation}
\label{ricci-soliton}
\big(r-R-\frac{\alpha'}{4}R^2\big)g_{ij}=\nabla_iX_j+\nabla_jX_i.
\end{equation}
If we moreover assume that \(X=-\nabla f\) for a function \(f(x,t)\), then we obtain
\begin{equation}
\label{gradient-ricci-soliton}
(R+\frac{\alpha'}{4}R^2-r)g_{ij}=2\nabla_i\nabla_j f.
\end{equation}
\begin{Dfn}
A metric \(g(t)\) on \(M\) is called 
second order Ricci soliton if it satisfies \eqref{ricci-soliton}.
Moreover, \(g(t)\) is called second order gradient Ricci soliton if it satisfies \eqref{gradient-ricci-soliton}.
\end{Dfn}
\begin{Bem}
It would be desirable to have an example of an explicit solution to \eqref{gradient-ricci-soliton}.
For the unnormalized version of the second order renormalization group flow
a generalization of the cigar solution of the Ricci flow was calculated in \cite{MR3202021}, Theorem 9.
Performing the same ansatz in our case gives an integro differential equation
and the existence of non-trivial solutions to this equation is more complicated.
\end{Bem}

Taking the trace of \eqref{gradient-ricci-soliton} yields
\begin{equation}
\label{curvature-potential}
R+\frac{\alpha'}{4}R^2-r=\Delta f.
\end{equation}
We call the function \(f\) the \emph{second order curvature potential}.
Note that \eqref{curvature-potential} is always solvable since the left hand side has vanishing integral.
We define the trace free part of \eqref{gradient-ricci-soliton} by
\begin{align}
M:=\nabla\nabla f-\frac{1}{2}g\Delta f.
\end{align}
It is easy to see that \(M=0\) is equivalent to the second order gradient Ricci soliton equation.

By a direct calculation it follows that for a gradient Ricci soliton we have
\begin{equation}
0=2\nabla^jM_{ij}=R\nabla_if+\big(1+\frac{\alpha'}{2}R\big)\nabla_i R
\end{equation}
and in the case that \(R>0\)
\[
\nabla(\log R+f+\frac{\alpha'}{2}R)=0.
\]
Ignoring the sign of the scalar curvature yields the conservation law
\begin{equation}
\frac{\alpha'}{4}R^2\nabla_if=\nabla_i\big(|\nabla f|^2+rf+R+\frac{\alpha'}{4}R^2\big)
\end{equation}
and combining both equations gives
\begin{equation}
0=\nabla_i\big(|\nabla f|^2+rf+R+\alpha'\frac{3}{8}R^2+\frac{\alpha'^2}{24}R^3\big).
\end{equation}

\begin{Bem}
In the case of the ``usual'' Ricci flow on closed surfaces,
gradient solitons on \(S^2\) correspond to metrics of constant curvature.
In the case of second order gradient Ricci solitons on \(S^2\) we have the following:
Taking the trace of \eqref{ricci-soliton} we get
\[
r-R-\frac{\alpha'}{4}R^2=\operatorname{div} X.
\]
Multiplying with \(R-r\) and integrating over \(S^2\) we find
\begin{align}
\label{second-order-kazdan}
\int_{S^2}\big(R+\frac{\alpha'}{4}R^2-r\big)(R-r\big)d\mu&=-\int_{S^2}\langle\nabla R,X\rangle d\mu=0,
\end{align}
where we applied the Kazdan-Warner identity (see \cite{MR0343205}) in the last step.
It is obvious, that metrics of constant curvature satisfy \eqref{second-order-kazdan}.
However, there may be additional solutions to \eqref{second-order-kazdan}.
\end{Bem}

We would now like to derive an evolution equation for the second order curvature potential \(f\).

\begin{Lem}
Let \((M,g_t)\) be a solution of the normalized second order renormalization group flow \eqref{second-order-flow}.
Then the following equation holds
\begin{align}
\label{evolution-laplace-partialt-f}
\Delta\frac{\partial f}{\partial t}=\big(1+\frac{\alpha'}{2}R\big)\Delta\Delta f+\frac{\alpha'}{4}R^2\Delta f+r\Delta f-\frac{\partial r}{\partial t}.
\end{align}
\end{Lem}
\begin{proof}
We differentiate \eqref{curvature-potential} with respect to \(t\) and find
\begin{align*}
\big(R+\frac{\alpha'}{4}R^2-r\big)\Delta f+\Delta\frac{\partial f}{\partial t}=(1+\frac{\alpha'}{2}R\big)\frac{\partial R}{\partial t}-\frac{\partial r}{\partial t}
=(1+\frac{\alpha'}{2}R\big)\big(\Delta\Delta f+R\Delta f\big)-\frac{\partial r}{\partial t},
\end{align*}
which proves the result.
\end{proof}
\begin{Bem}
In the case of the usual Ricci flow we can derive an evolution equation for the curvature potential \(f\).
Unfortunately, this does not seem to be possible here.
Interchanging derivatives gives
\[
\Delta\big(\frac{\partial f}{\partial t}-\big(1+\frac{\alpha'}{2}R\big)\Delta f-rf-\frac{\alpha'}{4}R^2f\big)=-\alpha'(\frac{1}{2}R\nabla R\nabla f+\frac{1}{4}f\Delta R^2+\nabla R\nabla\Delta f+\Delta R\Delta f)
-\frac{\partial r}{\partial t}.
\]
In general, the right hand side of this equation does not have a sign and we cannot extract any further information.
\end{Bem}

We do not get a nice evolution equation for \(f\).
Thus, we define another curvature potential:
\begin{Dfn}
We define the \emph{first order potential} \(w\) by
\begin{equation}
\label{def-first-order-potential}
\Delta w=R-a, \qquad a:=\frac{\int_MRd\mu}{\int_Md\mu}.
\end{equation}
\end{Dfn}
Note that \eqref{def-first-order-potential} is always solvable since the right hand side has vanishing integral.
Moreover, \eqref{def-first-order-potential} does not arise from a gradient soliton, which is different compared to the Ricci flow.

\begin{Lem}
Let \(w_0(x,t)\) be a first order potential for a solution \((M,g_t)\) of the normalized 
second order renormalization group flow \eqref{second-order-flow}.
Then there exists a function \(c(t)\) such that \(w:=w_0+c(t)\) satisfies the evolution equation
\begin{equation}
\label{evolution-w}
\frac{\partial w}{\partial t}=\big(1+\frac{\alpha'}{4}R+a\frac{\alpha'}{4}\big)\Delta w+af.
\end{equation}
\end{Lem}

\begin{proof}
We calculate
\begin{align*}
\frac{\partial}{\partial t}\Delta w&=(R+\frac{\alpha'}{4}R^2-r)\Delta w+\Delta\frac{\partial w}{\partial t}
=\frac{\partial R}{\partial t}
=\Delta(R+\frac{\alpha'}{4}R^2)+R(R+\frac{\alpha'}{4}R^2-r)
\end{align*}
and thus obtain
\[
\Delta\frac{\partial w}{\partial t}=\Delta\Delta w+\frac{\alpha'}{4}\Delta R^2+a(R+\frac{\alpha'}{4}R^2-r).
\]
Using
\(
\Delta R^2=\Delta(R(\Delta w+a)))=\Delta((R+a)\Delta w)
\)
we find
\[
\Delta\big(\frac{\partial w}{\partial t}-\big(1+\frac{\alpha'}{4}R+a\frac{\alpha'}{4}\big)\Delta w-af\big)=0.
\]
On a closed surface the only harmonic functions are constants, thus there is a function \(\gamma(t)\) satisfying
\[
\frac{\partial}{\partial t}w_0=\big(1+\frac{\alpha'}{4}R+a\frac{\alpha'}{4}\big)\Delta w_0+af+\gamma(t).
\]
The claim follows from setting \(c(t)=-\int_0^t\gamma(\tau)d\tau\)
and simplifying the right hand side.
\end{proof}

\begin{Bem}
Let us make some comments on the evolution equation \eqref{evolution-w}.
First of all, we do not have control over \(f\).
Secondly, it can easily be checked that \eqref{evolution-w} is parabolic for \(g\in \cM_+\).
Then, via the maximum principle, we have the estimate
\[
w\leq a\int_0^Tf(x,\tau)d\tau+w_0.
\]
\end{Bem}

We can also rewrite the equation for \(w\) as follows:
\begin{Bem}
Let \((M,g_t)\) be a solution of the normalized second order renormalization group flow \eqref{second-order-flow}.
Then the first order potential \(w\) satisfies
\[
\Delta\big(\frac{\partial w}{\partial t}-\big(1+\frac{\alpha'}{4}R+\frac{\alpha'}{4}a\big)\Delta w)-aw\big)=a\frac{\alpha'}{4}(R^2-\frac{1}{\vol(M,g_t)}\int R^2d\mu).
\]
\end{Bem}

\begin{Bem}
There holds the following relation between first and second order potential 
\begin{equation}
\Delta f=(1+\frac{\alpha'}{4}R)\Delta w+(1+\frac{\alpha'}{4}R)a-r.
\end{equation}
\end{Bem}

As a next step we calculate the evolution of \(|\nabla w|^2\) and \(|\nabla^2 w|^2\) under \eqref{second-order-flow}.
\begin{Lem}
Let \((M,g_t)\) be a solution of the normalized second order renormalization group flow \eqref{second-order-flow}.
The quantity \(|\nabla w|^2\) evolves by
\begin{equation}
\label{evolution-nabla-w}
\frac{\partial}{\partial t}|\nabla w|^2=(1+\frac{\alpha'}{2}R)\Delta|\nabla w|^2
-2(1+\frac{\alpha'}{2}R)|\nabla^2w|^2-(r+\frac{\alpha'}{4}R^2)|\nabla w|^2+2a\langle\nabla f,\nabla w\rangle.
\end{equation}
\end{Lem}
\begin{proof}
Using \eqref{evolution-w} we calculate
\begin{align*}
\frac{\partial}{\partial t}|\nabla w|^2=&\big(R+\frac{\alpha'}{4}R^2-r\big)|\nabla w|^2+2\langle\nabla\frac{\partial w}{\partial t},\nabla w\rangle \\
=&\big(R+\frac{\alpha'}{4}R^2-r\big)|\nabla w|^2+2\big(1+\frac{\alpha'}{4}R+a\frac{\alpha'}{4}\big)\langle\nabla\Delta w,\nabla w\rangle \\
&+\frac{\alpha'}{2}\nabla R\Delta w\nabla w+2a\langle\nabla f,\nabla w\rangle.
\end{align*}
Interchanging derivatives and using 
\begin{align*}
\frac{\alpha'}{2}\nabla R\Delta w\nabla w&=-\frac{\alpha'}{4}R(R-a)|\nabla w|^2+\frac{\alpha'}{2}(R-a)\langle\Delta\nabla w,\nabla w\rangle
\end{align*}
we obtain
\[
\frac{\partial}{\partial t}|\nabla w|^2=2(1+\frac{\alpha'}{2}R)\langle\Delta\nabla w,\nabla w\rangle-(r+\frac{\alpha'}{4}R^2)|\nabla w|^2
+2a\langle\nabla f,\nabla w\rangle,
\]
which gives the result.
\end{proof}

From \eqref{evolution-nabla-w} it becomes clear that we can only use the first order potential to
derive estimates in the case when \(a=0\) since we do not need any control over \(\nabla f\)
in that case.

\begin{Lem}
Let \((M,g_t)\) be a solution of the normalized second order renormalization group flow \eqref{second-order-flow}.
Then the the quantity \(|\nabla^2 w|^2\) evolves by
\begin{align}
\frac{\partial}{\partial t}|\nabla^2 w|^2=&(1+\frac{\alpha'}{2}R)\Delta|\nabla^2 w|^2-2(1+\frac{\alpha'}{2}R)|\nabla^3w|^2
-2(R+\frac{3}{4}\alpha'R^2+r)|\nabla^2w|^2 \\
\nonumber &+2R(1+\frac{\alpha'}{2}R)(R-a)^2+\alpha'\nabla R\nabla R\nabla^2w-2a\nabla^2f\nabla^2w.
\end{align}
\end{Lem}

\begin{proof}
We calculate
\begin{align*}
\frac{\partial}{\partial t}\nabla_i\nabla_jw=&\nabla_i\nabla_j\frac{\partial w}{\partial t}-\big(\frac{\partial}{\partial t}\Gamma^k_{ij}\big)\nabla_kw \\
=&\nabla_i\nabla_j\big((1+\frac{\alpha'}{4}R+\frac{\alpha'}{4}a)\Delta w-af\big)\\
&+\frac{1}{2}\big(1+\frac{\alpha'}{2}R\big)
\big(\nabla_iR\nabla_jw+\nabla_jR\nabla_iw-\langle\nabla R,\nabla w\rangle g_{ij}\big),
\end{align*}
where we applied \eqref{evolution-christoffel-second-order} for the evolution of the Christoffel symbols.
Using
\[
\nabla_i\nabla_j(R\Delta w)=2R\nabla_i\nabla_j\Delta w+2\nabla_jR\nabla_iR-a\nabla_i\nabla_j\Delta w
\]
and the identity
\begin{align*}
\nabla_i\nabla_j\Delta w=&\Delta\nabla_i\nabla_jw-\frac{1}{2}\big(\nabla_iR\nabla_jw+\nabla_jR\nabla_iw-\langle\nabla R,\nabla w\rangle g_{ij}\big) 
-2R(\nabla_i\nabla_jw-\frac{1}{2}(\Delta w)g_{ij})
\end{align*}
we obtain
\begin{align*}
\frac{\partial}{\partial t}\nabla_i\nabla_jw=&(1+\frac{\alpha'}{2}R)\big(\Delta\nabla_i\nabla_jw
-2R\nabla_i\nabla_jw+R(R-a)g_{ij} \big) 
+\frac{\alpha'}{2}\nabla_iR\nabla_jR-a\nabla_i\nabla_jf.
\end{align*}
The claim then follows from
\[
\frac{\partial}{\partial t}|\nabla^2w|^2=2(R+\frac{\alpha'}{4}R^2-r)|\nabla^2 w|^2+2\langle\frac{\partial}{\partial t}\nabla^2w,\nabla^2w\rangle.
\]
\end{proof}
\subsection{Surfaces with zero Euler Characteristic}
On a surface with \(\chi(M)=0\) the evolution equations \eqref{evolution-w} and \eqref{evolution-nabla-w} simplify:
\begin{align}
\nonumber\frac{\partial w}{\partial t}=&\big(1+\frac{\alpha'}{4}R\big)\Delta w, \\
\label{nabla-w-torus}\frac{\partial}{\partial t}|\nabla w|^2=&(1+\frac{\alpha'}{2}R)\Delta|\nabla w|^2
-2(1+\frac{\alpha'}{2}R)|\nabla^2w|^2-(r+\frac{\alpha'}{4}R^2)|\nabla w|^2.
\end{align}
Note that all terms involving \(f\) and \(\nabla f\) drop out. 
Thus, we may expect that we can use the curvature potential \(w\) in this case
to derive estimates. Moreover, since the Euler characteristic is zero an upper bound on \(R\)
should also yield a lower bound.

First of all, by applying the maximum principle to \eqref{nabla-w-torus} with \(g\in \cM_+\), we obtain
\begin{equation}
|\nabla w_T|^2+\int_0^T(1+\frac{\alpha'}{2}R)R^2dt\leq|\nabla w_0|^2.
\end{equation}
Moreover, rewriting \eqref{nabla-w-torus} as
\begin{align*}
\frac{\partial}{\partial t}\log|\nabla w|^2=&(1+\frac{\alpha'}{2}R)\Delta\log|\nabla w|^2+(1+\frac{\alpha'}{2}R)\big|\nabla\log|\nabla w|^2\big|^2
-2(1+\frac{\alpha'}{2}R)\frac{|\nabla^2w|^2}{|\nabla w|^2}\\
&-(r+\frac{\alpha'}{4}R^2)
\end{align*}
and by the maximum principle with \(g\in \cM_+\) we obtain
\begin{equation}
\log|\nabla w_T|^2+\int_0^T(r+\frac{\alpha'}{4}R^2)dt\leq\log|\nabla w_0|^2
\end{equation}
giving us control over the curvature in this case.

\section{Some Tools for surfaces with positive Euler Characteristic}
In the case of positive Euler characteristic \(\chi(M)>0\) we have
\[
r=\frac{2\pi\chi(M)}{\vol(M,g)}+\frac{1}{\vol(M,g)}\frac{\alpha'}{4}\int R^2d\mu>0.
\]
Since 
\[
-r|\nabla w|^2+2a\langle\nabla f,\nabla w\rangle\leq a|\nabla f|^2
\]
we find for the curvature potential in this case
\[
\frac{\partial}{\partial t}|\nabla w|^2\leq(1+\frac{\alpha'}{2}R)\Delta|\nabla w|^2
-2(1+\frac{\alpha'}{2}R)|\nabla^2w|^2-\frac{\alpha'}{4}R^2|\nabla w|^2+a|\nabla f|^2.
\]
Unfortunately, we again do not have control over \(\nabla f\) to obtain an estimate.

\subsection{Positive Scalar curvature}
Throughout this section we assume that \(R>0\).
Note that combining \eqref{evolution-dmu} and \eqref{evolution-R} gives
\begin{equation}
\label{evolution-measure}
\frac{\partial}{\partial t}(Rd\mu)=\Delta (R+\frac{\alpha'}{4}R^2)d\mu,
\end{equation}
which also means that \(\frac{\partial}{\partial t}\chi(M)=0\).

\begin{Dfn}
We define the entropy for the second order renormalization group flow as
\begin{equation}
N=\int_M \big(R\log R+\frac{\alpha'}{4}R^2\big)d\mu.
\end{equation}
\end{Dfn}

\begin{Lem}
Let \((M,g_t)\) be a solution of the normalized second order renormalization group flow \eqref{second-order-flow}.
Then the second order surface entropy evolves as
\begin{equation}
\frac{\partial N}{\partial t}=\int_M \big(-\big(\frac{\alpha'}{2}\sqrt{R}+\frac{1}{\sqrt{R}}\big)^2|\nabla R|^2+\big(R+\frac{\alpha'}{4}R^2-r\big)^2\big)d\mu.
\end{equation}
\end{Lem}
\begin{proof}
First of all, we calculate
\begin{align*}
\frac{\partial}{\partial t}\int_M R\log Rd\mu=&\int_M\big(\frac{\partial}{\partial t}\log R\big)Rd\mu+\int_M\log R\frac{\partial}{\partial t}(Rd\mu) \\
=&\int_MR\big(R+\frac{\alpha'}{4}R^2-r\big)d\mu+\int_M\log R\Delta (R+\frac{\alpha'}{4}R^2)d\mu
\end{align*}
using \eqref{evolution-measure}. Moreover, we find
\[
\frac{\partial}{\partial t}\frac{\alpha'}{4}\int R^2d\mu=-\frac{\alpha'}{2}\int(1+\frac{\alpha'}{2}R)|\nabla R|^2d\mu+\frac{\alpha'}{4}\int R^2\big(R+\frac{\alpha'}{4}R^2-r\big)d\mu
\]
\end{proof}
and adding up both contributions yields the result.

\begin{Prop}
Let \((M,g_t)\) be a solution of the normalized second order renormalization group flow \eqref{second-order-flow}.
Then the second order surface entropy satisfies
\begin{align}
\frac{\partial N}{\partial t}=-\int_M\frac{|\nabla R+R\nabla f+\frac{\alpha'}{2}R\nabla R|^2}{R}d\mu-2\int_M|M|^2d\mu\leq 0,
\end{align}
where \(M\) is the trace-free part of the Hessian of \(f\).
\end{Prop}
\begin{proof}
A direct calculation shows
\begin{align*}
\int_M\frac{|\nabla R+R\nabla f+\frac{\alpha'}{2}R\nabla R|^2}{R}d\mu=
\int_M\big(&|\nabla R|^2\big(\frac{\alpha'}{2}\sqrt{R}+\frac{1}{\sqrt{R}}\big)^2+R|\nabla f|^2 \\
&-2\big(R+\frac{\alpha'}{4}R^2\big)\big(R+\frac{\alpha'}{4}R^2-r\big)\big)d\mu.
\end{align*}
Moreover, we have (see the proof of Prop. 5.39 in \cite{MR2061425} for a derivation)
\[
-2\int_M|M|^2d\mu=\int_M\big(R|\nabla f|^2-\big(R+\frac{\alpha'}{4}R^2-r\big)^2\big)d\mu.
\]
Subtracting both equations, we have
\begin{align*}
-2\int_M|M|^2d\mu-&\int_M\big(\frac{|\nabla R+R\nabla f+\frac{\alpha'}{2}R\nabla R|^2}{R}\big)d\mu=\\
&-\int_M\big(|\nabla R|^2\big(\frac{\alpha'}{2}\sqrt{R}+\frac{1}{\sqrt{R}}\big)^2+\int_M\big(R+\frac{\alpha'}{4}R^2-r\big)^2d\mu.
\end{align*}
The claim then follows from the last Lemma.
\end{proof}
\begin{Cor}
Let \((M,g_t)\) be a solution of the normalized second order renormalization group flow \eqref{second-order-flow}
on a closed surface with \(R>0\). 
Then the entropy is strictly-decreasing unless \(g\) is a second order gradient soliton.
\end{Cor}
\begin{proof}
If \(\frac{dN}{dt}=0\) at a time \(t_0\in [0,T)\), then \(M(\cdot,t_0)=0\), giving the result.
\end{proof}

\begin{Bem}
The second order surface entropy satisfies
\begin{equation}
N(g)=r\vol(M)-2\pi\chi(M)+\int_MR\log Rd\mu.
\end{equation}
\end{Bem}

\begin{Bem}
If the curvature changes its sign we could try to 
consider the solution \(s\) of
\begin{align*}
\frac{d}{dt}s=s(s+\frac{\alpha'}{4}s^2-r)
\end{align*}
and set \(P:=R-s\).
We define the modified entropy as
\begin{align*}
\hat N(g_t,s_t):=\int_M\big(P\log P+\frac{\alpha'}{4}P^2\big),
\end{align*}
which evolves by
\begin{align}
\label{evolution-modified-entropy}
\frac{d}{dt}\hat N=&
\int_M\big(-\big(\frac{\alpha'}{2}\sqrt{P}+\frac{1}{\sqrt{P}}\big)^2|\nabla P|^2-\big(\frac{\alpha'}{2}\frac{s}{P}+\frac{\alpha'^2}{4}s\big)|\nabla P|^2\big)d\mu \\
\nonumber&+\int_MP(R+s-r+s\log P)d\mu+\frac{\alpha'}{4}\int_MP(R^2+s^2+Rs+\log P(s^2+Rs)d\mu \\
\nonumber&+\frac{\alpha'}{4}\int_M(P^2(R+2s-r))d\mu +\frac{\alpha'^2}{16}\int_MP^2(R^2+2s^2+2Rs)d\mu.
\end{align}
It becomes obvious from \eqref{evolution-modified-entropy} that we will not get a nice
estimate for the modified entropy.
\end{Bem}

\begin{Bem}
For the analysis of the Ricci flow on closed surfaces of positive Euler characteristic a useful quantity is 
\begin{equation*}
G:=R^2+t|\nabla R|^2.
\end{equation*}
Under the second order renormalization group flow the quantity \(G\) evolves as
\begin{align*}
\label{evolution-G}
\frac{\partial G}{\partial t}=&\big(1+\frac{\alpha'}{2}R\big)\Delta G-2\big(1+\frac{\alpha'}{2}R\big)t|\nabla^2R|^2+\alpha't\nabla|\nabla R|^2\nabla R \\
&\nonumber+|\nabla R|^2\big(-1+t(4R-3r+\frac{5}{4}\alpha'R^2+\alpha'\Delta R)\big)
+2R^2(R+\frac{\alpha'}{4}R^2-r). 
\end{align*}
and, unfortunately, we do not get a nice estimate by the maximum principle.
\end{Bem}

\begin{Bem}
In turns out that if we would have a time-dependent bound on the scalar curvature \(R\),
we could use the second order entropy to establish a uniform bound on the scalar curvature
similar to the methods used for the Ricci flow, see \cite{MR2061425}, p.141.
\end{Bem}

\subsection{Towards a Harnack inequality}
For simplicity we again only consider the case \(R>0\).
From the equation for second order gradient Ricci solitons \eqref{gradient-ricci-soliton} we obtain
\[
\nabla(\log R+f+\frac{\alpha'}{2}R)=0.
\]
Setting 
\begin{equation}
L:=\log R+\frac{\alpha'}{2}R 
\end{equation}
we may define
\begin{Dfn}
The second order differential Harnack quantity is defined by
\begin{equation}
Q:=\Delta L+R+\frac{\alpha'}{4}R^2-r.
\end{equation}
\end{Dfn}
Note that on a second order gradient soliton of positive curvature
the second order differential Harnack quantity \(Q\) satisfies \(Q=0\).
Hence, one should expect that \(Q\) satisfies a nice evolution equation.
\begin{Lem}
Let \((M,g_t)\) be a solution of the normalized second order renormalization group flow \eqref{second-order-flow}.
Then the quantity \(L\) evolves by
\begin{equation}
\frac{\partial L}{\partial t}=(1+\frac{\alpha'}{2}R)\Delta L+|\nabla L|^2+\big(1+\frac{\alpha'}{2}R\big)\big(R+\frac{\alpha'}{4}R^2-r\big).
\end{equation}

\begin{proof}
Using the evolution equation for the scalar curvature \eqref{evolution-R} we find
\[
\frac{\partial L}{\partial t}=(1+\frac{\alpha'}{2}R)\Delta L+(1+\alpha'R)|\nabla\log R|^2+\frac{\alpha'^2}{4}|\nabla R|^2
+\big(1+\frac{\alpha'}{2}R\big)\big(R+\frac{\alpha'}{4}R^2-r\big)
\]
and noting that
\[
|\nabla L|^2=(1+\alpha'R)\frac{|\nabla R|^2}{R^2}+\frac{\alpha'^2}{4}|\nabla R|^2
\]
yields the assertion.
\end{proof}

\end{Lem}
\begin{Cor}
The second order differential Harnack quantity thus satisfies
\begin{equation}
(1+\frac{\alpha'}{2}R)Q=\frac{\partial L}{\partial t}-|\nabla L|^2.
\end{equation}
\end{Cor}
\begin{Lem}
Let \((M,g_t)\) be a solution of the normalized second order renormalization group flow \eqref{second-order-flow}.
Then the quantity \(Q\) evolves by
\begin{align}
\frac{\partial Q}{\partial t}=\Delta(1+\frac{\alpha'}{2}R)Q+2\langle\nabla Q,\nabla L\rangle
\nonumber+2\big|\nabla\nabla L+\frac{1}{2}\big(R+\frac{\alpha'}{4}R^2-r\big)g\big|^2
+\big(\frac{\alpha'}{4}R^2+r\big)Q-\frac{\partial r}{\partial t}.
\end{align}
\end{Lem}
\begin{proof}
We calculate
\begin{align*}
\frac{\partial Q}{\partial t}=&\frac{\partial}{\partial t}(\Delta L)+\frac{\partial}{\partial t}\big(R+\frac{\alpha'}{4}R^2\big)-\frac{\partial r}{\partial t}\\
=&\big(R+\frac{\alpha'}{4}R^2-r\big)\Delta L+\Delta\frac{\partial L}{\partial t}+R\frac{\partial L}{\partial t} -\frac{\partial r}{\partial t}\\
=&\big(R+\frac{\alpha'}{4}R^2-r\big)\Delta L+\Delta(|\nabla L|^2+(1+\frac{\alpha'}{2}R)Q)\\
&+R\big((1+\frac{\alpha'}{2}R)\Delta L+|\nabla L|^2+\big(1+\frac{\alpha'}{2}R\big)\big(R+\frac{\alpha'}{4}R^2-r\big)\big)-\frac{\partial r}{\partial t}.
\end{align*}
Note that
\begin{align*}
\Delta|\nabla L|^2
=&2|\nabla\nabla L|^2+R|\nabla L|^2+2\langle\nabla Q,\nabla L\rangle-2(1+\frac{\alpha'}{2}R)\langle\nabla R,\nabla L\rangle \\
=&2|\nabla\nabla L|^2-R|\nabla L|^2+2\langle\nabla Q,\nabla L\rangle
\end{align*}
and thus we arrive at 
\begin{align*}
\frac{\partial Q}{\partial t}=&\Delta(1+\frac{\alpha'}{2}R)Q+2\langle\nabla Q,\nabla L\rangle \\
&+\big(2R+\frac{3}{4}\alpha'R^2-r\big)\Delta L+2|\nabla\nabla L|^2+\big(R+\frac{\alpha'}{2}R^2\big)\big(R+\frac{\alpha'}{4}R^2-r\big)-\frac{\partial r}{\partial t}\\
=&\Delta((1+\frac{\alpha'}{2}R)Q)+2\langle\nabla Q,\nabla L\rangle+2\big|\nabla\nabla L+\frac{1}{2}\big(R+\frac{\alpha'}{4}R^2-r\big)g\big|^2 \\
&+\big(\frac{\alpha'}{4}R^2+r\big)(\Delta L+R+\frac{\alpha'}{4}R^2-r)-\frac{\partial r}{\partial t}.
\end{align*}
Finally, using the definition of \(Q\), we get the result.
\end{proof}

\begin{Cor}
By the maximum principle, we may estimate \(Q\) as follows
\begin{equation}
\label{estimate-Q}
\frac{\partial Q}{\partial t}\geq Q^2+\big(\frac{\alpha'}{2}\Delta R+\frac{\alpha'}{4}R^2+r\big)Q-\frac{\partial r}{\partial t}.
\end{equation}
\end{Cor}

\begin{Bem}
The estimate \eqref{estimate-Q} is somehow unsatisfactory since we need control over \(\Delta R\)
to gain an estimate over \(R\).
\end{Bem}

\bibliographystyle{plain}
\bibliography{mybib}
\end{document}